# On Semi-Analytical Integration Specified for Mass Matrix of Finite Elements.

E Hanukah

Faculty of Mechanical Engineering, Technion – Israel Institute of Technology, Haifa 32000, Israel

Email: eliezerh@tx.technion.ac.il

**Abstract**

Spatial numerical integration is essential for finite element analysis. Currently, numerical integration schemes, mostly based on Gauss quadrature, are widely used. Herein, we present an alternative semi-analytical approach for mass matrix evaluation, resulting in an accurate, efficient and easy-to-implement integration rule. To this end, integrands of mass matrix entries get separated to two multiplicative parts. The first depends on natural coordinates while the second depends also on element parameters (mesh). Second part is approximated using polynomials and function evaluations at sampling (integration) points, allowing later analytical integration to precompute the weight matrixes. Resulting formulas possess typical form such that computational efficiency equivalence to traditional schemes is demonstrated, namely our n-point formula is computationally equivalent to the conventional use of n-point quadrature. Specifically, one and four-point semi-analytical formulas are explicitly derived for an eight node brick element. Our one-point rule describes exactly constant metric and density elements; four-point scheme correctly admits the linear metric density multiplication. Dramatic superiority in terms of accuracy is established based on coarse mesh study. Wherefore test sets are generated with the help of ABAQUS$^{TM}$, and then analyzed via variety of formulas. Averaged absolute errors, together with minimum and maximum absolute errors of standard and the new formulas are listed and comparisons are presented.

Authors of the report intend to achieve the next goals: ● Introduction of an alternative way to compute the mass matrix of solid finite elements, including inhomogeneous continua, consistent as well as inconsistent mass matrices. ● New easy-to-implement, efficient, formulas specified for widely-used eight-node brick element. ● Demonstration of accuracy advantage based on coarse mesh generated with commercial software.

## 1. Introduction

Evaluation of mass coefficients, which require integration in the element domain, is most commonly obtained by numerical integration (e.g. [1, 2]). Computational cost is proportional to the number of integration point, therefore its reduction is vital, yet, adopted quadrature should (i) permit convergence (ii) maintain convergence rate which would occur if exact integration were used ([3]pp. 164).

In addition, it has been shown that in many cases symbolic computations combined with code generation give rise to significant speedup of traditional quadrature use (e.g.[4, 5]). Herein we show how the accuracy of the n-point conventional integration rule can be dramatically increased by preprocessing additional integrand information, namely by specifying conventional integration formula to the particular problem of mass matrix calculation. Our approach yields efficient, problem specific formulas. Efficiency follows from integration point reduction.

The idea of intelligently using symbolic, closed-form or semi-analytical integration combined with algebraic manipulations, to produce accurate and efficient matrices integration for particular finite elements isn't new (e.g. [6, 7]). For instance, stiffness matrixes has been successfully treated for a number of cases: plane elements [8-11], 3D bricks [9, 12, 13], triangles [14, 15] and 3D elastic and electro-elastic tetrahedral [15-20] etc.

Furthermore, hierarchical semi-analytical, displacement based approach, for which complete order polynomials are used, allows to model three dimensional finite bodies e.g. [21-23] resulting in new analytical solutions (e.g.[24]). Effective yet not a symbolic approach has been developed to speedup numerical quadrature for p-version FE matrices e.g. [25].

The gap addressed here is the mass matrices integration appropriate for solid elements with *varying metric* i.e. coarse and fine mesh elements allowed. Moreover, varying initial density (inhomogeneity) is also systematically treated. Our earlier findings [26-29] are systematized and generalized.

The term "semi-analytical" is adopted to indicate two key ideas: *analytical* integration is used to precompute the weights matrices; however, the resulting formulas are *not exact* (semi) due to prior polynomial approximation. Specifically, we separate the integrand to two multiplicative parts. First part is purely polynomial function with respect to natural coordinates. Second part is the metric (Jacobian determinant) density multiplication, which in addition to polynomial dependence on local coordinates, depends also on element characteristics (parameters) – nodal positions (mesh). Polynomial approximation to the second part is suggested yielding a series of ansatz functions multiplied by the function evaluations at certain sampling point. Then analytical integration is applied to precompute the weight matrices. Entries of the weight matrices are real numbers, similar to traditional weights, though, they have been specified for a particular problem – mass matrix integration.

Efficiency of the resulting formula is discussed based on comparison to conventional quadrature use. It is shown, that roughly the same amount of operations has to be performed to use traditional and the proposed formulas. In particular, the standard scheme requires slightly more multiplication while the new one requires slightly more weights to be remembered.

As a special example, a consistent mass matrix of standard eight node brick element is considered. Two formulas are offered in a ready for implementation manner: first requires one integration point and assumes Constant Metric*Density (CMD) in the domain; second is a four-point scheme which assumes Linear Metric*Density (LMD). Importantly, both schemes are exact for homogeneous density and arbitrary parallelepiped mesh, since then the metric is constant. Thus, only coarse mesh accuracy study is necessary.

Comparative accuracy examination is carried out with the help of commercial software ABAQUS. Five curved bodies are meshed with hexahedrons; nodal positions together with connectivity arrays are adopted and processed independently. Along with our one and four-point rules, performance of conventional one, four and six-point schemes are recorded. In terms of average error, our CMD (one-point) over-performs even the six point scheme, while our LMD formula demonstrates very accurate and robust behavior. Average accuracy for poor meshes remains below 1% error. Details of minimum and maximum errors are recorded and establish four-point (LMD) formula as appropriate for all practical mesh (fine and coarse).

The outline of the note as follows. Section 2 introduces the semi-analytical approach to mass matrix evaluation. In Section 3, the method is applied to a particular problem – eight node brick element. One (CMD) and four-point (LMD) formulas are explicitly developed and presented. Section 4 contains numerical accuracy study. To this end, hexahedral coarse mesh is generated with the help of ABAQUS. Performance comparison between conventional schemes and the new ones is conducted. Section 6 summarizes and records our conclusions. Appendix contains eight-node solid brick details: shape functions, metric etc.

## 2. Method.

Standard isoparametric formulation (e.g.[30] pp.104) of solid (mass conserving) element, yield a consistent symmetric mass matrix [M], which entries are denoted by $M_{ij} = M_{ji}$

$$M_{ij} = \int_v \rho N_i N_j dv = \int_V \rho_0 N_i N_j dV = \int_V N_i N_j (\rho_0 J) d\square \ , \ d\square = d\xi d\eta d\zeta \ , \ (i,j=1,..,n_{nodes}) \quad (1)$$

Where $n_{nodes}$ stands for a number of nodes in the element, $\{\xi, \eta, \zeta\}$ symbolize natural coordinates, $v, V$ and $\rho, \rho_0$ denote actual and initial element volumes and densities accordingly. Shape functions $N_i, (i=1,..,n_{nodes})$ merely polynomials with respect to natural coordinates e.g.(11), moreover mesh independent. Inhomogeneous initial configuration i.e. varying initial density $\rho_0 = \rho_0(\xi, \eta, \zeta)$ is considered. Metric (or Jacobian) - determinant of the Jacobi matrix of global-local coordinate's transformation is denoted by $J$. Importantly, $J$ is a polynomial functions with regard to natural coordinates in addition a mesh dependent term $J = J(\xi, \eta, \zeta; X_{mi})$ (e.g.(12)), where $X_{mi}$ $(m=1,2,3, i=1,..,n_{nodes})$ indicate components of initial nodal positions $\mathbf{X}_i$ $(i=i,..,n_{nodes})$ in terms of the global Cartesian coordinate system $\{\mathbf{e}_1, \mathbf{e}_2, \mathbf{e}_3\}$, i.e.

$\mathbf{X}_i = X_{mi}\mathbf{e}_m$ (m = 1, 2, 3). Here and throughout the text, bold symbols stand for vector quantities; traditional summation convention on repeated indexes is implied. Finally it is noted, that the forthcoming discussion is equally applicable to inconsistent mass matrices such as diagonal (lumped) or mixed (e.g. [31]).

In (1) we sharply distinguish between two multiplicative parts of the integrand; first part depends purely on natural coordinates $N_i N_j = N_i(\xi,\eta,\zeta) N_j(\xi,\eta,\zeta)$, $(i, j = 1,..,n_{nodes})$ while the other is a mesh dependent term $\rho_0 J = \rho_0(\xi,\eta,\zeta) J(\xi,\eta,\zeta; X_{mi})$. Contrary to the metric $J$, which explicitly depends on mesh $X_{mi}$ (12), initial density $\rho_0$ is a continuum property; yet, we lack convenient density function so it is appropriate to treat it together with the metric. In addition, it is emphasized that polynomials such as shape functions or their multiplications can be easily analytically integrated. Therefore we propose to replace mesh dependent function $\rho_0 J$ by polynomial approximation $\rho_0 J \approx \hat{N}_p(\xi,\eta,\zeta)(\rho_0 J)|_{\text{point } p}$. Where $(\rho_0 J)|_{\text{point } p}$ indicates function evaluation at a specific sampling (integration) point. Clearly, the more sampling points one uses, the higher approximation order and accuracy of the approximation. Substitution to (1) result in $M_{ij} \approx (\int_V N_i N_j \hat{N}_p d\square)(\rho_0 J)|_{\text{point } p}$ thus

$$M_{ij} \approx \sum_{p=1}^{n_p} \hat{M}_{ijp} \rho_{0p} J_p$$

$$\hat{M}_{ijp} = \int_V N_i N_j \hat{N}_p d\square \ , \ J_p = J|_{\text{point } p} \ , \ \rho_{0p} = \rho_0|_{\text{point } p} \ , \ (i, j = 1,..,n_{nodes})$$

(2)

Where $n_p$ represents the number of integration points. Essentially, entries of symmetric weight matrices $\hat{M}_{ijk} = \hat{M}_{jik}$ are precomputed numbers similar to standard weights in conventional rules. Implementation of $\hat{M}_{ijk}$ in code is also similar, numeric values apply for all elements of the same kind. It is vital to draw the connection to standard way of using quadrature, such that efficiency equivalence will be clear. Currently, entries (1) are estimated as

$$M_{ij} \approx \sum_{p=1}^{n_p} w_p (N_i N_j \rho_0 J)|_{\text{point } p} =$$

$$\sum_{p=1}^{n_p} w_p N_i(\xi_p,\eta_p,\zeta_p) N_j(\xi_p,\eta_p,\zeta_p) \rho_0(\xi_p,\eta_p,\zeta_p) J(\xi_p,\eta_p,\zeta_p; X_{mk})$$

$$(m = 1,2,3, i, j = 1,..,n^{nodes})$$

Using shorter notation, standard numerical integration (e.g. [32]pp.75) takes the form

$$M_{ij} \approx \sum_{p=1}^{n_p} M_{ijp} \rho_{0p} J_p \quad , \quad M_{ijp} = w_p N_{ip} N_{jp} \quad , \quad (m=1,2,3, i,j=1,..,n_{nodes}) \quad (3)$$

$$N_{ip} = N_i(\xi_p, \eta_p, \zeta_p) \quad , \quad \rho_{0p} = \rho_0(\xi_p, \eta_p, \zeta_p) \quad , \quad J_p = J(\xi_p, \eta_p, \zeta_p; X_{mk})$$

Obviously, for a general case (3) and (2) are not equal, nevertheless, forms are identical, hence roughly the same amount of computations have to be performed for equal $n_p$. In particular, to evaluate $M_{563}$ subroutine takes $w_3, N_{53}, N_{63}$ from the memory and multiplies them (3), while $\hat{M}_{563}$ is just pooled out from the memory (2). On the other hand, evaluation of $M_{ijp}$ requires only $w_p (p=1,..,n_p)$ and $N_{ip} (i=1,..,n_{nodes}, p=1,..,n_p)$ to be stored in the memory, consequently $n_p + n_{nodes}*n_p$ real numbers, while for symmetric $\hat{M}_{ijp} = \hat{M}_{jip}$ $(i,j=1,..,n_{nodes}, p=1,..,n_p)$ one has to store $n_p*(n_{nodes}*n_{nodes}+n_{nodes})/2$ numbers. Importantly, this additional memory demand is absolutely negligible in modern hardware. Moreover, as it is demonstrated in (9),(10) most of these numbers are zero or repeat one another.

### 3. Specific example – eight node brick element.

Details of (2) semi-analytical integration formula - sample points and weight matrices, have to be explicitly stated for each formula individually. Herein, we consider widely used hexahedral element (e.g.[30]pp.119) providing all necessary details for implementation. Evaluation points (e.g.[2]pp.230) for one and four-point schemes

$$1-\text{point}: \quad (\xi_1, \eta_1, \zeta_1) = (0,0,0) \quad (4)$$

$$4-\text{point}: \quad (\xi_1, \eta_1, \zeta_1) = (0,0,0) \quad , \quad (\xi_2, \eta_2, \zeta_2) = (\frac{1}{10},0,0)$$
$$(\xi_3, \eta_3, \zeta_3) = (0,\frac{1}{10},0) \quad , \quad (\xi_4, \eta_4, \zeta_4) = (0,0,\frac{1}{10}) \quad (5)$$

Perhaps, the simplest (zero order) model, for the metric density multiplication $\rho_0 J \approx \hat{N}_p \rho_{0p} J_p$ $(p=1,..,n_p)$ is the Constant Metric Density (CMD) $\rho_0 J \approx \rho_{01} J_1$ for which only one point and one ansatz functions $\hat{N}_p (p=1,..,n_p)$ is needed

$$M_{ij}^{CMD} = \hat{M}_{ij1} \rho_{01} J_1 \quad , \quad \hat{M}_{ij1} = \int_{-1}^{+1} \int_{-1}^{+1} \int_{-1}^{+1} N_i N_j \hat{N}_1^{CMD} d\xi d\eta d\zeta \quad , \quad \hat{N}_1^{CMD} = 1 \quad , \quad (i,j=1,..,8) \quad (6)$$

With the help of the above in common with shape functions (11)

$$\hat{M}_{ijl} = \frac{1}{27}\begin{pmatrix} 8 & 4 & 2 & 4 & 4 & 2 & 1 & 2 \\ 4 & 8 & 4 & 2 & 2 & 4 & 2 & 1 \\ 2 & 4 & 8 & 4 & 1 & 2 & 4 & 2 \\ 4 & 2 & 4 & 8 & 2 & 1 & 2 & 4 \\ 4 & 2 & 1 & 2 & 8 & 4 & 2 & 4 \\ 2 & 4 & 2 & 1 & 4 & 8 & 4 & 2 \\ 1 & 2 & 4 & 2 & 2 & 4 & 8 & 4 \\ 2 & 1 & 2 & 4 & 4 & 2 & 4 & 8 \end{pmatrix} \quad (7)$$

Linear Metric Density (LMD) model in 3D space requires four sampling points besides, linear shape functions

$$M_{ij}^{LMD} = M_{ijp}\rho_{0p}J_p \;,\; (p,k=1,..,4) \;,\; \hat{N}_p^{LMD}\Big|_{\text{point } k} = \delta_{pk}$$
$$N_1^{LMD} = 1 - 10\xi - 10\eta - 10\zeta \;,\; N_2^{LMD} = 10\xi \;,\; N_3^{LMD} = 10\eta \;,\; N_4^{LMD} = 10\zeta \quad (8)$$

$\delta_{pk}$ denotes Kronecker symbol. Linear ansatz function $N_p^{LMD}$ together with analytical integration produce

$$\hat{M}_{ij1} = \frac{1}{27}\begin{pmatrix} 128 & 44 & 12 & 44 & 44 & 12 & 1 & 12 \\ 44 & 48 & 4 & 12 & 12 & 4 & -8 & 1 \\ 12 & 4 & -32 & 4 & 1 & -8 & -36 & -8 \\ 44 & 12 & 4 & 48 & 12 & 1 & -8 & 4 \\ 44 & 12 & 1 & 12 & 48 & 4 & -8 & 4 \\ 12 & 4 & -8 & 1 & 4 & -32 & -36 & -8 \\ 1 & -8 & -36 & -8 & -8 & -36 & -112 & -36 \\ 12 & 1 & -8 & 4 & 4 & -8 & -36 & -32 \end{pmatrix} ,\; \hat{M}_{ij2} = \frac{1}{27}\begin{pmatrix} -40 & 0 & 0 & -20 & -20 & 0 & 0 & -10 \\ 0 & 40 & 20 & 0 & 0 & 20 & 10 & 0 \\ 0 & 20 & 40 & 0 & 0 & 10 & 20 & 0 \\ -20 & 0 & 0 & -40 & -10 & 0 & 0 & -20 \\ -20 & 0 & 0 & -10 & -40 & 0 & 0 & -20 \\ 0 & 20 & 10 & 0 & 0 & 40 & 20 & 0 \\ 0 & 10 & 20 & 0 & 0 & 20 & 40 & 0 \\ -10 & 0 & 0 & -20 & -20 & 0 & 0 & -40 \end{pmatrix} \quad (9)$$

$$\hat{M}_{ij3} = \frac{1}{27}\begin{pmatrix} -40 & -20 & 0 & 0 & -20 & -10 & 0 & 0 \\ -20 & -40 & 0 & 0 & -10 & -20 & 0 & 0 \\ 0 & 0 & 40 & 20 & 0 & 0 & 20 & 10 \\ 0 & 0 & 20 & 40 & 0 & 0 & 10 & 20 \\ -20 & -10 & 0 & 0 & -40 & -20 & 0 & 0 \\ -10 & -20 & 0 & 0 & -20 & -40 & 0 & 0 \\ 0 & 0 & 20 & 10 & 0 & 0 & 40 & 20 \\ 0 & 0 & 10 & 20 & 0 & 0 & 20 & 40 \end{pmatrix} ,\; \hat{M}_{ij4} = \frac{1}{27}\begin{pmatrix} -40 & -20 & -10 & -20 & 0 & 0 & 0 & 0 \\ -20 & -40 & -20 & -10 & 0 & 0 & 0 & 0 \\ -10 & -20 & -40 & -20 & 0 & 0 & 0 & 0 \\ -20 & -10 & -20 & -40 & 0 & 0 & 0 & 0 \\ 0 & 0 & 0 & 0 & 40 & 20 & 10 & 20 \\ 0 & 0 & 0 & 0 & 20 & 40 & 20 & 10 \\ 0 & 0 & 0 & 0 & 10 & 20 & 40 & 20 \\ 0 & 0 & 0 & 0 & 20 & 10 & 20 & 40 \end{pmatrix} \quad (10)$$

Although, more accurate models can be suggested by invoking more polynomials and integration points, still, based on the next section, LMD is fairly accurate for all practical purpose. Furthermore, 4 sampling points proposed here for LMD, do not match the standard 4 point. According to our findings, the proposed points will on average, in the nonlinear least square sense (e.g. [33, 34]) result in slightly lower errors. Nevertheless, several different sets of 4 point, including the traditional one (e.g.[30]pp.122, [2]pp.229) has been checked out, resulting in insignificant loss of accuracy. Linear metric density model uses complete order polynomials that create accurate scheme, for which, point locations of minor importance.

## 4. Comparative accuracy study

In the previous section, two integration formulas for a consistent mass matrix $M_{ij}$ $(i, j = 1,..,n_{nodes})$ of an eight node brick element have been detailed: one-point CMD and four-point LMD. It is important to perform comparative accuracy study based on coarse mesh. Fine mesh study is redundant, as it was shown in section 2; our CMD formula is exact for constant metric density i.e. for homogeneous initial density (in the element domain) our one-point (CMD) scheme is exact for all parallelepiped shapes of initial mesh, regardless of aspect ratio, as all parallelepipeds (even very skewed ones in all 3 dimensions) have constant metric. Our four-point (LMD) is exact for constant and linear metric density multiplication.

Commercial FE software ABAQUS has been used to produce five test sets of coarse hexahedral mesh, based on which, quantitative accuracy comparison of different schemes is conducted. For all five cases, homogeneous initial density configuration $\rho_0 = \text{const}$ has been assumed. We deliberately generated solid bodies whose mesh, generated by built-in algorithm, will be of poor quality, as we intended to check the scheme's performance for the worst "case scenario". To illustrate mesh coarseness, metric $J$ has been evaluated at all the elements (at hundreds of points inside element domains) and it was found that in the first case metrics drops below zero at 3 elements, for case 2 at 2 elements and for case 4 at 7 elements. Cases 3 and 5 do not have elements with negative metric (which is clearly reflected in higher accuracy). See Figures: 1–4 for visual details; number of nodes and elements apiece also recorded.

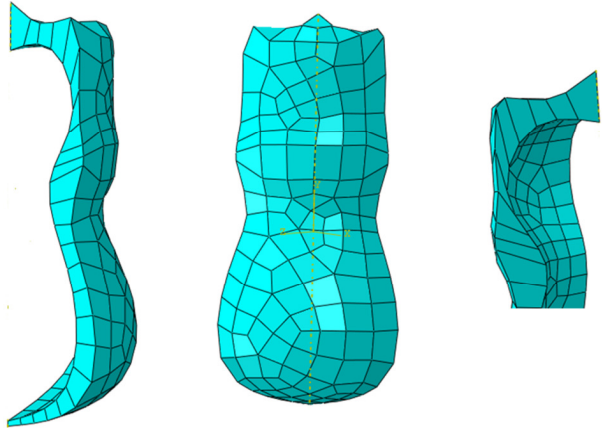

Figure 1: Three views of the first case. Consist of 328 nodes and 139 elements.

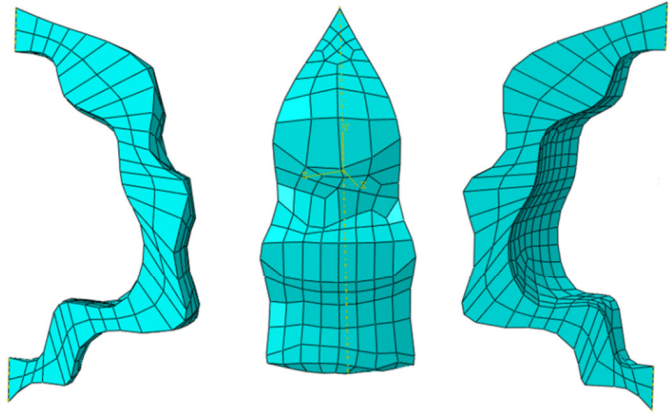

Figure 2: Three views of the second case. Contain of 594 nodes and 330 elements.

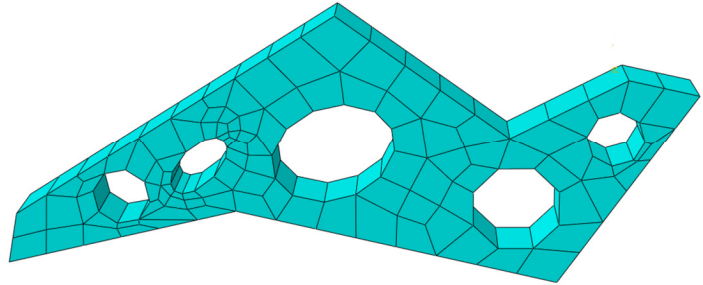

Figure 3: Third case. Involve 282 nodes and 100 elements.

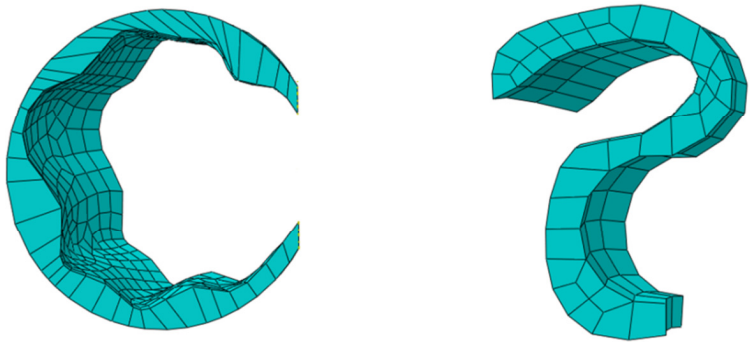

Figure 4: The 4$^{th}$ the 5$^{th}$ cases. Contain 666 nodes with 292 elements and 244 nodes with 111 elements respectively.

Nodal positions and connectivity array has been adopted from ABAQUS .inp file. All later analysis has been conducted using MAPLE$^{TM}$ hand written code. First of all, an exact mass matrix components computation was performed and used as a reference for error evaluation. Error of an individual matrix entry was calculated as an absolute value of $100*(estimated - exact)/exact$. Then, a list of an average absolute error was prepared, where each value is an average evaluated among element mass matrix components. Consequently, every element was indicated by a number – average absolute error. As a result, for the five test cases, average, maximum and minimum errors (average errors) have been recorded in the below Tables1 and 2.

Our Semi-Analytical (SA) CMD formula is a one-point scheme, hence, conventional one point rule has been examined for comparison purpose. Similarly, traditional four point scheme has been considered as a standard equivalence of our four-point LMD rule. In addition, performance of two different six-point formulas is included to once again illuminate accuracy advantages of the proposed schemes. First 6-point formula is suggested by FE textbook ([30] pp. 122) while the other (denoted by star) can be found in the quadrature "Bible" ([2] pp.229).

| CASE # | Elem tot # | Average absolute error % | | | | | |
|---|---|---|---|---|---|---|---|
| | | Standard quadrature rules | | | | New SA rules | |
| | | 1-p | 4-p | 6-p | 6*-p | CMD-1p | LMD-4p |
| 1 | 139 | 69.2 | 27.7 | 24.3 | 10.9 | 9.2 | 0.9 |
| 2 | 330 | 69.0 | 27.4 | 24.1 | 10.8 | 8.1 | 0.8 |
| 3 | 100 | 69.1 | 27.2 | 24.2 | 10.8 | 7.4 | 0.00000 |
| 4 | 292 | 69.4 | 27.7 | 24.4 | 10.9 | 9.2 | 0.8 |
| 5 | 111 | 68.9 | 27.3 | 24.0 | 10.8 | 6.4 | 0.4 |

Table 1: Average absolute error of conventional 1,4,6-point schemes together with semi-analytical 1 and 4-point rules, for all the five mesh cases.

Following Table 1, CMD-1point formula over-performs in terms of an average error, the standard 1, 4 and 6-point quadrature, while our LMD-4point scheme is less than 1% average error for all the cases. Especially we pay attention to the $3^{rd}$ and $5^{th}$ cases, for which even though the mesh is disported, no singular elements were found. Performance of our rules demonstrated the best accuracy in those cases. Moreover, Maximum errors (see Table 2) were also the lowest in those cases. Cases 1,2,4 serve as a "worst case scenario" since an experienced user will not tolerate such a poor mesh. Nevertheless, our SA 4-point scheme is fairly accurate based on the above average error and on the maximum-minimum table below.

| CASE # | Elem tot # | MAX-MIN error range % | | | | | | | | | | |
|---|---|---|---|---|---|---|---|---|---|---|---|---|
| | | Standard quadrature rules | | | | | | | | New SA rules | | |
| | | 1-p | | 4-p | | 6-p | | 6*-p | | CMD-1p | | LMD-4p |
| 1 | 139 | 75.6 | 68.1 | 31.7 | 27.0 | 28.9 | 23.4 | 12.5 | 10.5 | 26.0 | 1.9 | 4.6 | 0.1 |
| 2 | 330 | 75.1 | 68.3 | 30.9 | 27.0 | 31.6 | 23.5 | 12.3 | 10.1 | 29.4 | 0.8 | 10.6 | 0.1 |
| 3 | 100 | 74.6 | 68.4 | 28.8 | 27.1 | 27.2 | 23.6 | 11.9 | 10.4 | 19.3 | 0.3 | 0.0 | 0.0 |
| 4 | 292 | 97.1 | 68.3 | 37.4 | 27.0 | 37.5 | 23.5 | 14.3 | 10.1 | 56.4 | 0.9 | 5.1 | 0.0 |
| 5 | 111 | 72.9 | 68.3 | 28.2 | 27.1 | 26.4 | 23.6 | 11.5 | 10.6 | 19.4 | 0.7 | 1.1 | 0.0 |

Table 2: Maximum and minimum absolute error of conventional 1,4,6-point schemes together with semi-analytical 1 and 4-point rules, for all the five mesh cases.

Minimal errors of CMD formula are above zero which reflects the fact that there were no parallelepiped shaped elements of any kind; otherwise, our CMD scheme would be exact. In addition, the suggested LMD rule is fairly accurate at all times. Comparison of semi-analytical to equivalent conventional quadrature use reveals overwhelming superiority of the suggested method (see Table 1, 2). Importantly, our four point formula, demonstrate low errors even for coarse mesh.

## 5. Summary and Conclusions

Numerical integration of mass matrix of solid finite elements has been considered. Conventional quadrature use is challenged by the proposed semi-analytical method. Following the advocated approach, integrand is spited to two multiplicative parts; the first is polynomial with respect to natural coordinates, while the second depends explicitly on nodal positions. Polynomial approximation for the second part, which is the metric density multiplication, has to be suggested, and then analytical integration is used to calculate the weights matrices. Weight matrices entries are mesh independent, precomputed real numbers, similar to traditional weights. Resulting integration schemes take familiar easy to implement form, with major difference being the weighting. Weight matrices have to be specified for the particular element and number of integration points. Computational equivalence of semi-analytical n-point rule to traditional n-point rule is discussed.

In particular, consistent mass matrix, of widely used, eight node brick element has been studied. One and four-point rules were explicitly stated in a ready to be implemented manner i.e. integration points and weight matrices are detailed. Our one-point (Constant Metric Density - CMD) formula is exact for constant density elements with arbitrary parallelepiped configuration, while our four-point scheme admits linear metric density (LMD). Accuracy of the offered rules has been examined; five sets of coarse mesh were generated using ABAQUS. Our one-point scheme over-performs, in terms of average accuracy, traditional 1, 4 and 6-point schemes, while the average error of our 4-point formula remains below 1% even for very poor mesh. Examination of maximum and minimum errors supports once again substantial accuracy advantage.

**Appendix**

Shape functions of a standard eight node hexahedral element are given by

$$\begin{aligned}
N_1 &= (1-\zeta-\eta+\eta\zeta-\xi+\xi\zeta+\xi\eta-\xi\eta\zeta)/8 \;,\; N_2 = (1-\zeta-\eta+\eta\zeta+\xi-\xi\zeta-\xi\eta+\xi\eta\zeta)/8 \\
N_3 &= (1-\zeta+\eta-\eta\zeta+\xi-\xi\zeta+\xi\eta-\xi\eta\zeta)/8 \;,\; N_4 = (1-\zeta+\eta-\eta\zeta-\xi+\xi\zeta-\xi\eta+\xi\eta\zeta)/8 \\
N_5 &= (1+\zeta-\eta-\eta\zeta-\xi-\xi\zeta+\xi\eta+\xi\eta\zeta)/8 \;,\; N_6 = (1+\zeta-\eta-\eta\zeta+\xi+\xi\zeta-\xi\eta-\xi\eta\zeta)/8 \\
N_7 &= (1+\zeta+\eta+\eta\zeta+\xi+\xi\zeta+\xi\eta+\xi\eta\zeta)/8 \;,\; N_8 = (1+\zeta+\eta+\eta\zeta-\xi-\xi\zeta-\xi\eta-\xi\eta\zeta)/8
\end{aligned} \quad (11)$$

Material point X occupies location $\mathbf{X}$ inside the element domain $-1 \leq \xi, \eta, \zeta \leq +1$, $\mathbf{X}$ is represented by $\mathbf{X} = N_i \mathbf{X}_i$ $(i=1,..,8)$, where $\mathbf{X}_i$ are nodal positions (mesh) $\mathbf{X}_i = X_{mi} \mathbf{e}_m$ $(m=1,2,3, i=1,..,n_{nodes})$. Clearly, the above representation of $\mathbf{X}$ is mesh dependent. Jacobian determinant (metric) of global-local coordinate's transformation

$$J = \mathbf{X},_\xi \times \mathbf{X},_\eta \cdot \mathbf{X},_\zeta = \begin{vmatrix} (\mathbf{X} \cdot \mathbf{e}_1),_\xi & (\mathbf{X} \cdot \mathbf{e}_1),_\eta & (\mathbf{X} \cdot \mathbf{e}_1),_\zeta \\ (\mathbf{X} \cdot \mathbf{e}_2),_\xi & (\mathbf{X} \cdot \mathbf{e}_2),_\eta & (\mathbf{X} \cdot \mathbf{e}_2),_\zeta \\ (\mathbf{X} \cdot \mathbf{e}_3),_\xi & (\mathbf{X} \cdot \mathbf{e}_3),_\eta & (\mathbf{X} \cdot \mathbf{e}_3),_\zeta \end{vmatrix} > 0$$

Where comma denotes partial derivative; $|\bullet|$ stand for determinant operator. Following the above, together with algebraic manipulations, the metric of an eight node brick admits the next variable separated form

$$\begin{aligned}
J = \tilde{J}_0 &+ \\
&\xi \tilde{J}_1 + \eta \tilde{J}_2 + \zeta \tilde{J}_3 + \\
&\xi\eta \tilde{J}_4 + \xi\zeta \tilde{J}_5 + \eta\zeta \tilde{J}_6 + \xi^2 \tilde{J}_7 + \eta^2 \tilde{J}_8 + \zeta^3 \tilde{J}_9 + \\
&\xi\eta\zeta \tilde{J}_{10} + \xi^2\eta \tilde{J}_{11} + \xi\eta^2 \tilde{J}_{12} + \xi^2\zeta \tilde{J}_{13} + \eta^2\zeta \tilde{J}_{14} + \xi\zeta^2 \tilde{J}_{15} + \eta\zeta^2 \tilde{J}_{16} + \\
&\xi^2\eta\zeta \tilde{J}_{17} + \xi\eta^2\zeta \tilde{J}_{18} + \xi\eta\zeta^2 \tilde{J}_{19} \\
\tilde{J}_k &= \tilde{J}_k(X_{mi}) \;,\; (k=0,..,19, m=1,2,3, i=1,..,8)
\end{aligned} \quad (12)$$

Evidently, metric is a fourth order function with respect to natural coordinates; constant metric, linear and quadratic terms include complete order polynomials, whereas the third and fourth orders are incomplete. Likewise, an exact model of the metric using sampling points should include the above monomials together with 20 evaluation points.